\documentclass[10pt]{amsart}
\usepackage[dvipdfm]{graphicx,color}
\usepackage{url}
\usepackage{mathrsfs}
\newtheorem{theorem}{Theorem}[section]
\newtheorem{lemma}[theorem]{Lemma}

\theoremstyle{definition}
\newtheorem{definition}[theorem]{Definition}

\theoremstyle{remark}

\numberwithin{equation}{section}

\title[A Rellich Type Theorem for Discrete Schr{\"o}dinger Operators]{A Rellich Type Theorem for Discrete Schr{\"o}dinger Operators}
\author{Hiroshi ISOZAKI}
\address{Division of Mathematics,
University of Tsukuba
\\
Tsukuba, 305-8571, JAPAN 
\\ isozakih@math.tsukuba.ac.jp, hmorioka@math.tsukuba.ac.jp}
\author{Hisashi Morioka}
\subjclass[2000]{Primary~81U40, Secondary~47A40}
\keywords{Schr\"{o}dinger operator, square lattice, Rellich's theorem}

\begin{document}
\baselineskip 14pt
\maketitle

\begin{abstract}
An analogue of Rellich's theorem is proved for discrete Laplacian on square lattice, and applied to show unique continuation property on certain domains as well as  non-existence of embedded eigenvalues for discrete Schr{\"o}dinger operators.
\end{abstract}


\section{Introduction}

The Rellich type theorem for the Helmholtz equation is the following assertion (\cite{Re43}): Suppose $u\in H_{loc}^2 ( {\bf R}^d )$ satisfies
 \begin{equation*}
(-\Delta -\lambda )u =0  \quad \text{in} \quad \{ |x|>R_0 \} ,
\end{equation*} 
 for some constants $\lambda $, $R_0 >0 $, and
\begin{equation*}
 u(x) = o(|x|^{-(d-1)/2}), \quad |x|\to\infty.
\end{equation*}
Then $u(x)=0 $ on $\{ |x| >R_0 \}$.

This theorem has been extended to a broad class of Schr{\"o}dinger operators, since it implies the non-existence of eigenvalues embedded in the continuous spectrum (see e.g. \cite{Ka59}, \cite{Ro69}, \cite{Ag70}),  and also  plays an important role in the proof of limiting absorption principle which yields the absolute continuity of the continuous spectrum (see e.g. \cite{Ei62}, \cite{IkSa72}). The Rellich type theorem states a local property at infinity of solutions. Namely, it  proves  $u(x) = 0$ on $\{|x|>R_1\}$ for some $R_1 > R_0$.  By the unique continuation property, it then follows that $u(x) = 0$ for $|x| > R_0$.  In the theory of linear partial differential equations (PDE), the Rellich type theorem can be regarded as the problem of division in the momentum space. In fact, given a linear PDE with constant coefficients $P(D)u = f$, $f$ being compactly supported, the Fourier transform leads to the algebraic equation $P(\xi)\widetilde u(\xi) = \widetilde f(\xi)$, where $\widetilde u(\xi)$ denotes the Fourier transform of $u(x)$. If $P(\xi)$ divides $\widetilde f(\xi)$, $u$ is compactly supported due to the Paley-Wiener theorem. This approach was pursued by Treves \cite{Tre60}, and then developed by Littman \cite{Lit66}, \cite{Lit70}, H{\"o}rmander \cite{Hor70} and Murata \cite{Mur76}. One should note that Besov spaces appear naturally through these works.
In this paper, we shall consider its extension to the discrete case.

Throughout the paper, we shall assume that $d \geq 2$. Let ${\bf Z}^d = \{n = (n_1,\cdots,n_d)$ $ ; n_i \in {\bf Z}\}$ be the square lattice, and $e_{1} = (1,0,\cdots,0)$, $\cdots$, $e_{d} = (0,\cdots,0,1)$ the standard
bases of ${\bf Z}^d$. 
The discrete Laplacian $\Delta_{disc}$ is defined by
\begin{equation}
 \big(\Delta_{disc}\widehat u\big)(n) =  \frac{1}{4}\sum_{j=1}^{d} \big( {\widehat u}(n + e_{j}) +
{\widehat u}(n - e_{j}) \Big) - \frac{d}{2}\,\widehat u(n)
\nonumber
\end{equation}
for a sequence $\{\widehat u(n)\}_{n\in{\bf Z}^d}$.
Our main theorem is the following.


\begin{theorem}\label{rellich}
Let $\lambda \in (0,d)$ and $R_0>0$. Suppose that a sequence $\{\widehat u(n)\}$, defined for $\{ n\in {\bf Z}^d \ ; \ |n| \geq R_0 \} $, satisfies
\begin{equation}
(- \Delta_{disc} - \lambda)\widehat u = 0 \quad  \text{in} \quad \{ n\in {\bf Z}^d \ ; \ |n| > R_0 \} ,
\label{S1Equation}
\end{equation}
\begin{equation}
\lim_{R\to\infty}\frac{1}{R}\sum_{R_0<|n|<R}|\widehat u(n)|^2 = 0.
\label{S1DecayCond}
\end{equation}
Then there exists $R_1 > R_0$ such that $\widehat u(n) = 0$ for $|n| > R_1$.
\end{theorem}

Note that the spectrum of $ \widehat{H}_0 = -\Delta_{disc} $ is equal to $[0,d] $ and it is absolutely continuous (see e.g. \cite{IsKo}).

A precursor of this theorem is given in the proof of Theorem 9 of Shaban-Vainberg \cite{Sha}. Their purpose is to compute the asymptotic expansion of the resolvent $\widehat{R}_0 (\lambda \pm i0 )=( -\Delta_{disc} -\lambda \mp i0 )^{-1}$ on ${\bf Z}^d$ and to find the associated radiation condition which implies the uniqueness of the solution to the discrete Helmholtz equation. 
 Let
\begin{equation}
{\bf T}^d = {\bf R}^d/(2\pi{\bf Z)}^d \cong  [-\pi,\pi]^d
\label{S1Torus}
\end{equation}
be the $d$-dimensional flat torus 
and put
\begin{equation}
h(x)=\frac{1}{2}\Big(d - \sum_{j=1}^d\cos x_j\Big) , \quad
M_{\lambda}=\big\{x\in {\bf T}^d \ ; \  h(x)= \lambda \big\}.
\label{h(x)Mlambda}
\end{equation}
$M_{\lambda }$ is called the \textit{Fermi surface} of the discrete Laplacian.
For $\lambda \in (0,d)\setminus {\bf Z}$, the Fermi surface $M_{\lambda}$ is a $(d-1)$-dimensional smooth submanifold of ${\bf T}^d$.
On the other hand, if $\lambda \in (0,d)\cap {\bf Z}$, $M_{\lambda }$ has some isolated singularities.
Passing to the Fourier series, $-\Delta_{disc} $ is unitary equivalent to the operator of multiplication by $h(x) $ on ${\bf T}^d$. 
Therefore, the computation of the behavior of  $\widehat R_0(\lambda\pm i0)$ boils down to that for  an integral  on $M_{\lambda}$.
For a compactly supported function $\widehat{f} \in \ell^2 ({\bf Z}^d )$ and 
$\lambda \in  (0,d)\setminus {\bf Z}$,  
the stationary phase method gives 
the following asymptotic expansion as $|n| \to \infty$:
\begin{gather} 
\begin{split}
&\big( \widehat{R}_0 (\lambda \pm i0 )\widehat{f} \big)(n) \\
&= |n|^{-(d-1)/2} \sum_{j} e^{\pm i n \cdot x_{\infty} ^{(j)} (\lambda , \omega_n )} a^{(j)}_{\pm} (\lambda , \omega_n ) +O(|n|^{-(d+1)/2}),
\end{split}
\label{SV_asymptotic}
\end{gather}
where $\omega_n = n/|n|$ is assumed to be \textit{non-singular} i.e. the Gaussian curvature does not vanish on all stationary phase points $x^{(j)}_{\infty} (\lambda,\omega_n) \in M_{\lambda}$ at which the normal of $M_{\lambda}$ is parallel to $\omega_n$.
It is then natural to define the radiation condition by using the first term of the above asymptotic expansion (\ref{SV_asymptotic}). 
To show the uniqueness of solutions to the discrete Helmholtz equation satisfying the radiation condition, they proved the assertion (which is buried in the proof actually): 

\medskip
\noindent
{\bf (S-V)} \ \textit{Let $\lambda \in (0,d)\setminus {\bf Z}$.
The solution of $(-\Delta_{disc} -\lambda )\widehat{u} =0 $ in $\{|n| >R_0 \}$ for $R_0 >0$, satisfying $\widehat{u} (n)=O(|n| ^{-(d+1)/2} )$, vanishes on $\{ |n| >R_1 \}$ for sufficiently large $R_1 >0$.}

\medskip

The new ingredient in the present paper is the following fact to be proved in \S 4.2:
Consider the equation
\begin{equation}
(h(x)-\lambda )u(x) = f(x), \quad {\rm on} \quad {\bf T}^d.
\label{S1EquationonTd}
\end{equation}
If the Fourier coefficients $\widehat u(n)$ of the distribution $u$ satisfy (\ref{S1DecayCond}) and $\widehat f(n)$ is compactly supported, then $u$ is smooth on ${\bf T}^d $ except for some null sets (in fact, the exceptional set is the set of discrete points), hence $f(x) = 0$ on $M_{\lambda }$, since $f$ is an analytic function. Our proof does not depend on the asymptotic expansion of the resolvent, and uses the optimal decay assumption 
(\ref{S1DecayCond}).
Moreover, it allows us to extend the theorem for interior threshold energy $\lambda \in (0,d)\cap {\bf Z}$, for which the Fermi surface has some singularities.
We derive some basic facts of the Fermi surface in \S 4.1.

Once we establish this fact, we can follow the arguments for proving the assertion (S-V) with some modifications to show that $\widehat u(n)$ is compactly supported. 
For the sake of completeness, in \S 4.2, we will also reproduce the proof of this part, which makes use of  basic facts in theories of functions of several complex variables and  algebraic geometry.

As applications of Theorem \ref{rellich}, we show in \S 2 non-existence of eigenvalues embedded in the continuous spectrum for $- \Delta_{disc}  + \widehat V$ in the whole space as well as in exterior domains. 
We also state the unique continuation property for exterior domains.

The result of the present paper is used as a key step in  \cite{IsMo} on the inverse scattering from the scattering matrix of a fixed energy for discrete Schr{\"o}dinger operators with compactly supported potentials.  In \cite{IsKo},   the inverse scattering  from all energies was studied by using complex Born approximation (see also \cite{Es}).

Function theory of  several complex variables and algebraic geometry have  already been utilized as powerful tools not only in linear PDE but also in the study of spectral properties for discrete Schr{\"o}dinger operators or periodic problems. See e.g.   Eskina \cite{Es},  Kuchment-Vainberg \cite{KuVa}, G{\'e}rard-Nier \cite{GeNi98}.

We give some remarks about notation in this paper.
For $x,y \in {\bf R}^d $, $x\cdot y= x_1 y_1 +\cdots +x_d y_d $ denotes the ordinary scalar product, and $|x| = (x \cdot x)^{1/2}$ is the Euclidean norm.
Note that even for $n = (n_1,\cdots,n_d) \in {\bf Z}^d$, we use $|n| = (\sum_{i=1}^d|n_i|^2)^{1/2}$.
For a Banach spaces $X$, $ \mathbf{B} (X)$ denotes the totality of bounded operators on $X$.
For a self-adjoint operator $A$ on a Hilbert space, $\sigma (A)$, $\sigma _{ess} (A) $, $\sigma_{ac} (A)$ and $\sigma_p (A)$ denote  the spectrum, the essential spectrum, the absolutely continuous spectrum and the point spectrum of $A$, respectively. 
For a set $S$, $\,^{\#}S$ denotes the number of elements in $S$. We use the notation
\begin{equation}
\langle t\rangle = (1 + |t|^2)^{1/2}, \quad t \in {\bf R},
\nonumber
\end{equation}
For a positive real number $a $ and $z' \in {\bf C}$, we denote $ z \equiv z' \ (\mathrm{mod} \ a )$ as a representative element of the equivalence class $ \{ z \in {\bf C} \ ; \ z=z' +a N \ \text{for} \ N \in {\bf Z} \}$.

\subsection{Acknowledgement}
The authors are indebted to Evgeny Korotyaev for useful discussions and  encouragements. The second author is supported by the Japan Society for the Promotion of Science under the Grant-in-Aid for Research Fellow (DC2) No. 23110.

\section{Some applications of Theorem \ref{rellich}}

\subsection{Absence of embedded eigenvalues in the whole space}
Granting  Theorem 1.1, we state its applications in this section.
 The  Schr{\"o}dinger operator
${\widehat H}$ on ${\bf Z}^d$ is defined by
\begin{equation}
{\widehat H} = - \Delta_{disc} + {\widehat V},
\nonumber
\end{equation}
where $\widehat V$ is the multiplication operator :
\begin{equation}
({\widehat V}\,{\widehat u})(n) = \widehat V(n)\,{\widehat u}(n).
\nonumber
\end{equation}


\begin{theorem}\label{NonExistenceWholeSpace}
If $\widehat V(n) \in {\bf R}$ for all $n$, and there exists $R_0 > 0$ such that 
$\widehat V(n) = 0$ for $|n| > R_0$, we have
$\sigma_p(\widehat H)\cap (0,d) = \emptyset.$
\end{theorem}

Since $\widehat{V}$ is compactly supported, Weyl's theorem yields $\sigma_{ess}(\widehat H) = [0,d]$ (see \cite{IsKo}). Therefore, Theorem \ref{NonExistenceWholeSpace} asserts the non-existence of eigenvalues embedded in the continuous spectrum except for the endpoints of $\sigma_{ess} (\widehat{H})$.
We also note that \cite{HSSS} has given an example of embedded eigenvalue at the endpoint of $\sigma_{ess} (\widehat{H} )$.

\medskip
{\it Proof of Theorem \ref{NonExistenceWholeSpace}}.  Suppose that $\lambda \in (0,d)  \cap \sigma_p (\widehat{H} )$ and $\widehat{u} \in \ell^2 ({\bf Z}^d )$ is the associated  eigenfunction i.e. $(- \Delta_{disc} +\widehat{V}-\lambda )\widehat{u}=0$.
Putting $\widehat{f}= - \widehat V\widehat{u} $, which is compactly supported by the assumption of $\widehat V$, we have the equation 
\begin{equation*}
(- \Delta_{disc}-\lambda )\widehat{u}=\widehat{f} \quad \text{on} \quad {\bf Z}^d. 
\end{equation*}
Since $\widehat u \in \ell^2({\bf Z}^d)$, the condition (\ref{S1DecayCond}) is satisfied. Theorem \ref{rellich} then implies that $\widehat{u} $ is compactly supported.
Therefore, there exists $m_1 \in {\bf Z}$ such that
 $\widehat{u}(n)=0 $ if $n_1 \geq m_1$.
Using the equation $(- \Delta_{disc} +\widehat{V} -\lambda )\widehat{u} =0$, which holds on the whole ${\bf Z}^d$, we then have  
\begin{gather*} 
\begin{split}
&\frac{1}{4} \widehat{u} ( m_1 -1 ,n') \\
&=\left( -\Delta_{disc}^{(d-1)} +\widehat{V} (m_1 , n') -\lambda\right)\widehat{u}(m_1 ,n')+\frac{1}{2} \widehat{u} (m_1 , n') -\frac{1}{4} \widehat{u} (m_1 +1 ,n') \\
&=0 , 
\end{split} 
\end{gather*}
where $n'=(n_2 , \cdots , n_d )$ and $\Delta_{disc}^{(d-1 )} $ is the discrete Laplacian on ${\bf Z}^{d-1}$. Repeating this procedure, we have $\widehat u(n) = 0$ for all $n$, and completes the proof. \qed


\subsection{Unique continuation property}
The next problem we address is the unique continuation property. We begin with the explanation of the exterior problem. 
 A subset $\Omega \subset {\bf Z}^d$ is said to be {\it connected} if for any $m, n \in \Omega$, there exists a sequence $n^{(0)}, \cdots, n^{(k)} \in \Omega$ with $n^{(0)}=m, n^{(k)}=n$ such that for all $0 \leq \ell \leq k-1$, $|n^{(\ell)}-n^{(\ell+1)}|=1$. 
For a connected subset $\Omega \subset {\bf Z}^d$, we put
\begin{equation}
{\rm deg}_{\Omega}(n) =\, ^{\# }\{m \in \Omega\, ; \, |m-n|=1\}, \quad n \in \Omega.
\label{S1Degree}
\end{equation}
The interior $\stackrel{\circ}\Omega$ and the boundary $\partial \Omega$ are defined by
\begin{gather}
\stackrel{\circ}\Omega\, = \{n \in \Omega\, ; \, {\rm deg}_{\Omega}(n) = 2d\},
\label{S1InteriorOmega} \\
\partial \Omega = \{n \in \Omega\, ; \, {\rm deg}_{\Omega}(n) < 2d\}.
\label{S1PartialD}
\end{gather}
The normal derivative on $\partial \Omega $ is defined by 
\begin{equation}
\partial _{\nu} \widehat{u}(n)= \frac{1}{4} \sum_{ m\in \stackrel{\circ}\Omega , |m-n|=1 } \big( \widehat{u}(n)-\widehat{u} (m) \big) , \quad n\in \partial \Omega .
\label{S2normald}
\end{equation}
Then, for a bounded  connected subset $\Omega $, the following Green formula holds:
\begin{gather}
\begin{split}
&\sum_{ n\in \stackrel{\circ}\Omega } \Big( \big( \Delta_{disc} \widehat{u} \big) (n) \cdot \widehat{v} (n) - \widehat{u} (n)\cdot \big( \Delta_{disc} \widehat{v} \big) (n)\Big) \\
&= \sum_{ n\in \partial \Omega} \Big( (\partial _{\nu }\widehat{u} \big) (n) \cdot \widehat{v} (n) - \widehat{u} (n) \cdot \big( \partial _{\nu} \widehat{v} \big) (n) \Big).
\end{split}
\label{greenformula}
\end{gather}
Indeed,  the standard definition of  Laplacian on graph  is (see e.g. \cite{Du84})
\begin{equation}
-\big( \Delta_{disc}^{\Omega} \widehat{u} \big)(n):=
\left\{
\begin{split}
&- \big( \Delta_{disc} \widehat{u} \big)(n) , \quad  ( n\in \stackrel{\circ}\Omega ), \\
&\big( \partial _{\nu} \widehat{u}\big) (n),  \quad  ( n\in \partial \Omega) ,
\end{split}
\right.
\label{graph_laplacian}
\end{equation}
which yields
\begin{equation}
\sum_{n\in \Omega } \big( \Delta_{disc}^{\Omega} \widehat{u} \big)(n) \cdot \widehat{v}(n) = \sum_{n\in \Omega}\widehat{u}(n) \cdot \big( \Delta_{disc}^{\Omega} \widehat{v} \big) (n) , \quad \widehat{u} , \ \widehat{v} \in \ell^2 (\Omega ).
\label{S2symmetry}
\end{equation}
Splitting the sum (\ref{S2symmetry}) into two parts, the ones over $\stackrel{\circ}\Omega $ and over $\partial \Omega $, we have (\ref{greenformula}).

Let $\Omega_{ext}$ be a exterior domain, which means that there is a bounded set $\Omega_{int}$ such that $\Omega_{ext} = {\bf Z}^d \, \setminus   \stackrel{\circ}\Omega_{int}$.
We assume that $\Omega_{ext} $ is connected.
We consider the Schr{\"o}dinger operator 
\begin{equation}
\widehat H_{ext} = - \Delta_{disc} + \widehat V
\label{S1Hext}
\end{equation}
without imposing the boundary condition, where $\widehat V$ is a real-valued compactly supported potential.

Now suppose there exists a $\lambda \in (0,d)$, and $\widehat u$ satisfying  (\ref{S1DecayCond}) and
\begin{equation}
(\widehat H_{ext} - \lambda)\widehat u = 0 \quad {\rm in}\quad \stackrel{\circ}\Omega_{ext}.
\label{S1ExtScroedingerEq}
\end{equation}
 By Theorem \ref{rellich}, $\widehat u$ vanishes near infinity. However, in the discrete case  the unique continuation property of Laplacian does not hold in general. It depends on the shape of the domain. To guarantee it, we introduce the following {\it cone condition}. For $1 \leq i \leq d$ and $n \in {\bf Z}^d$, let $C_{i,\pm}(n)$ be the cone defined by
\begin{equation}
C_{i,\pm}(n) = \Big\{ m \in {\bf Z}^d\, ; \, \sum_{k\neq i}|m_k-n_k| \leq \pm (m_i - n_i) \Big\}.
\label{S1Conecond}
\end{equation}


\begin{definition}\label{ConeCond}
An exterior domain $\Omega_{ext}$ is said to satisfy a {\it cone condition} if for any $n \in\Omega_{ext}$, there is a cone $C_{i,+}(n)$ (or $C_{i,-}(n)$) such that 
$C_{i,+}(n) \subset \Omega_{ext}$ (or $C_{i,-} (n)\subset \Omega_{ext}$).
\end{definition}

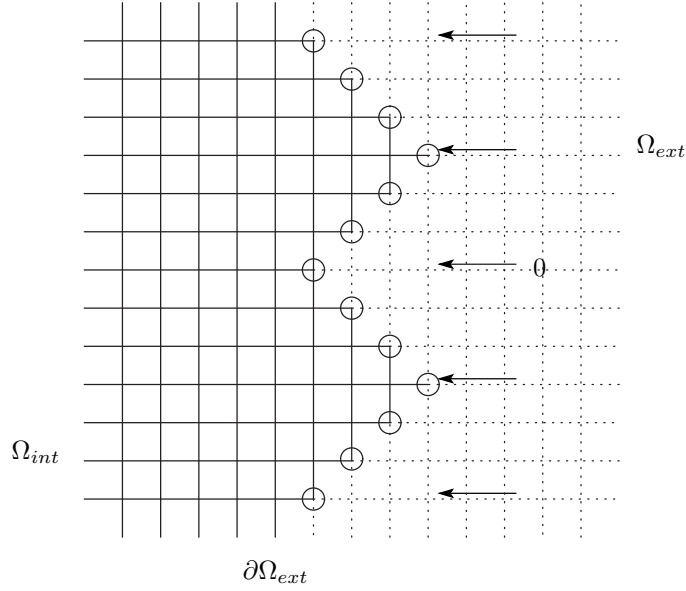
\begin{figure}[t]
\centering
\unitlength 0.1in
\begin{picture}( 32.7000, 29.1000)(  6.2000,-31.1000)
%
{\color[named]{Black}{%
\special{pn 8}%
\special{pa 2800 1000}%
\special{pa 1000 1000}%
\special{fp}%
\special{pa 2600 1200}%
\special{pa 1000 1200}%
\special{fp}%
\special{pa 1000 1400}%
\special{pa 2400 1400}%
\special{fp}%
\special{pa 2200 1600}%
\special{pa 1000 1600}%
\special{fp}%
\special{pa 1000 1800}%
\special{pa 2400 1800}%
\special{fp}%
\special{pa 1000 2000}%
\special{pa 2600 2000}%
\special{fp}%
\special{pa 1000 2200}%
\special{pa 2800 2200}%
\special{fp}%
\special{pa 2600 2400}%
\special{pa 1000 2400}%
\special{fp}%
\special{pa 1000 2600}%
\special{pa 2400 2600}%
\special{fp}%
\special{pa 2200 2800}%
\special{pa 1000 2800}%
\special{fp}%
\special{pa 1000 800}%
\special{pa 2600 800}%
\special{fp}%
\special{pa 2400 600}%
\special{pa 1000 600}%
\special{fp}%
\special{pa 1000 400}%
\special{pa 2200 400}%
\special{fp}%
\special{pa 2600 800}%
\special{pa 2600 1200}%
\special{fp}%
\special{pa 2400 600}%
\special{pa 2400 1400}%
\special{fp}%
\special{pa 2200 1600}%
\special{pa 2200 400}%
\special{fp}%
\special{pa 2000 200}%
\special{pa 2000 3000}%
\special{fp}%
\special{pa 2200 2800}%
\special{pa 2200 1600}%
\special{fp}%
\special{pa 2400 1800}%
\special{pa 2400 2600}%
\special{fp}%
\special{pa 2600 2400}%
\special{pa 2600 2000}%
\special{fp}%
\special{pa 1800 3000}%
\special{pa 1800 200}%
\special{fp}%
\special{pa 1600 200}%
\special{pa 1600 3000}%
\special{fp}%
\special{pa 1400 3000}%
\special{pa 1400 200}%
\special{fp}%
\special{pa 1200 200}%
\special{pa 1200 3000}%
\special{fp}%
}}%
%
{\color[named]{Black}{%
\special{pn 8}%
\special{ar 2200 400 58 58  0.0000000 6.2831853}%
}}%
%
{\color[named]{Black}{%
\special{pn 8}%
\special{ar 2400 600 58 58  0.0000000 6.2831853}%
}}%
%
{\color[named]{Black}{%
\special{pn 8}%
\special{ar 2600 800 58 58  0.0000000 6.2831853}%
}}%
%
{\color[named]{Black}{%
\special{pn 8}%
\special{ar 2800 1000 58 58  0.0000000 6.2831853}%
}}%
%
{\color[named]{Black}{%
\special{pn 8}%
\special{ar 2600 1200 58 58  0.0000000 6.2831853}%
}}%
%
{\color[named]{Black}{%
\special{pn 8}%
\special{ar 2400 1400 58 58  0.0000000 6.2831853}%
}}%
%
{\color[named]{Black}{%
\special{pn 8}%
\special{ar 2200 1600 58 58  0.0000000 6.2831853}%
}}%
%
{\color[named]{Black}{%
\special{pn 8}%
\special{ar 2400 1800 58 58  0.0000000 6.2831853}%
}}%
%
{\color[named]{Black}{%
\special{pn 8}%
\special{ar 2600 2000 58 58  0.0000000 6.2831853}%
}}%
%
{\color[named]{Black}{%
\special{pn 8}%
\special{ar 2800 2200 58 58  0.0000000 6.2831853}%
}}%
%
{\color[named]{Black}{%
\special{pn 8}%
\special{ar 2600 2400 58 58  0.0000000 6.2831853}%
}}%
%
{\color[named]{Black}{%
\special{pn 8}%
\special{ar 2400 2590 58 58  0.0000000 6.2831853}%
}}%
%
{\color[named]{Black}{%
\special{pn 8}%
\special{ar 2200 2800 58 58  0.0000000 6.2831853}%
}}%
\put(33.5000,-16.3000){\makebox(0,0)[lb]{$0$}}%
%
{\color[named]{Black}{%
\special{pn 8}%
\special{pa 3260 1570}%
\special{pa 2860 1570}%
\special{fp}%
\special{sh 1}%
\special{pa 2860 1570}%
\special{pa 2928 1590}%
\special{pa 2914 1570}%
\special{pa 2928 1550}%
\special{pa 2860 1570}%
\special{fp}%
\special{pa 3260 970}%
\special{pa 2860 970}%
\special{fp}%
\special{sh 1}%
\special{pa 2860 970}%
\special{pa 2928 990}%
\special{pa 2914 970}%
\special{pa 2928 950}%
\special{pa 2860 970}%
\special{fp}%
\special{pa 3260 2170}%
\special{pa 2860 2170}%
\special{fp}%
\special{sh 1}%
\special{pa 2860 2170}%
\special{pa 2928 2190}%
\special{pa 2914 2170}%
\special{pa 2928 2150}%
\special{pa 2860 2170}%
\special{fp}%
\special{pa 3260 2770}%
\special{pa 2860 2770}%
\special{fp}%
\special{sh 1}%
\special{pa 2860 2770}%
\special{pa 2928 2790}%
\special{pa 2914 2770}%
\special{pa 2928 2750}%
\special{pa 2860 2770}%
\special{fp}%
\special{pa 3260 370}%
\special{pa 2860 370}%
\special{fp}%
\special{sh 1}%
\special{pa 2860 370}%
\special{pa 2928 390}%
\special{pa 2914 370}%
\special{pa 2928 350}%
\special{pa 2860 370}%
\special{fp}%
}}%
\put(6.2000,-26.1000){\makebox(0,0)[lb]{$\Omega_{int}$}}%
\put(18.3000,-32.4000){\makebox(0,0)[lb]{$\partial \Omega_{ext}$}}%
%
{\color[named]{Black}{%
\special{pn 8}%
\special{pa 2200 200}%
\special{pa 2200 400}%
\special{dt 0.045}%
\special{pa 2200 400}%
\special{pa 3800 400}%
\special{dt 0.045}%
\special{pa 2400 600}%
\special{pa 3800 600}%
\special{dt 0.045}%
\special{pa 3800 800}%
\special{pa 2600 800}%
\special{dt 0.045}%
\special{pa 2800 1000}%
\special{pa 3800 1000}%
\special{dt 0.045}%
\special{pa 2600 1200}%
\special{pa 3800 1200}%
\special{dt 0.045}%
\special{pa 3800 1400}%
\special{pa 2400 1400}%
\special{dt 0.045}%
\special{pa 2200 1600}%
\special{pa 3800 1600}%
\special{dt 0.045}%
\special{pa 3800 1800}%
\special{pa 2400 1800}%
\special{dt 0.045}%
\special{pa 2600 2000}%
\special{pa 3800 2000}%
\special{dt 0.045}%
\special{pa 3800 2200}%
\special{pa 2800 2200}%
\special{dt 0.045}%
\special{pa 2600 2400}%
\special{pa 3800 2400}%
\special{dt 0.045}%
\special{pa 3800 2600}%
\special{pa 2400 2600}%
\special{dt 0.045}%
\special{pa 2200 2800}%
\special{pa 3800 2800}%
\special{dt 0.045}%
\special{pa 2200 2800}%
\special{pa 2200 3000}%
\special{dt 0.045}%
\special{pa 2400 3000}%
\special{pa 2400 2600}%
\special{dt 0.045}%
\special{pa 2600 3000}%
\special{pa 2600 2400}%
\special{dt 0.045}%
\special{pa 2800 3000}%
\special{pa 2800 2200}%
\special{dt 0.045}%
\special{pa 3000 3000}%
\special{pa 3000 200}%
\special{dt 0.045}%
\special{pa 2400 200}%
\special{pa 2400 600}%
\special{dt 0.045}%
\special{pa 2600 200}%
\special{pa 2600 800}%
\special{dt 0.045}%
\special{pa 2800 200}%
\special{pa 2800 1000}%
\special{dt 0.045}%
\special{pa 2800 1000}%
\special{pa 2800 2200}%
\special{dt 0.045}%
\special{pa 2600 1200}%
\special{pa 2600 2000}%
\special{dt 0.045}%
\special{pa 2400 1400}%
\special{pa 2400 1800}%
\special{dt 0.045}%
\special{pa 3200 3000}%
\special{pa 3200 200}%
\special{dt 0.045}%
\special{pa 3400 200}%
\special{pa 3400 3000}%
\special{dt 0.045}%
\special{pa 3600 3000}%
\special{pa 3600 200}%
\special{dt 0.045}%
}}%
\put(38.9000,-10.1000){\makebox(0,0)[lb]{$\Omega_{ext}$}}%
\end{picture}%
\caption{The zig zag type boundary.}
\label{Fig1}
\end{figure}

Examples of the domain satisfying this cone condition are
\begin{itemize}
\item $\left(\Omega_{ext}\right)^c$ = a rectangular polyhedron = $\{n \in {\bf Z}^d\, ; \, |n_i| \leq a_i, \ i = 1, \cdots, d\}$,
\item $\left(\Omega_{ext}\right)^c$ = a rhombus = $\{n \in {\bf Z}^d\, ; \, \sum_{i=1}^d|n_i| \leq C\}$,
\item a domain with zigzag type boundary (see Figure \ref{Fig1}).
\end{itemize}

\begin{theorem} \label{UniqueConti}
Let $\widehat H_{ext}$ be a Schr{\"o}dinger operator in an exterior domain $\Omega_{ext} \subset {\bf Z}^d$ with compactly supported potential. Suppose $\Omega_{ext}$ satisfies the cone condition. If there exist $\lambda \in (0,d) $ and $\widehat u$ satisfying (\ref{S1ExtScroedingerEq}) and (\ref{S1DecayCond}), then $\widehat u = 0$ on $\Omega_{ext}$.
\end{theorem}

Proof. Take  any $n \in \Omega_{ext}$. 
By the cone condition, there is a cone, say $C_{1,+}(n)$, such that $C_{1,+}(n) \subset \Omega_{ext}$. 
By Theorem \ref{rellich}, there is $k_1$ such that 
$\widehat u(m) = 0$ for $m \in C_{1,+}(n)$, $k_1 < m_1$. 
Arguing as in the proof of Theorem \ref{NonExistenceWholeSpace}, we have $\widehat u(k_1,m') = 0$, $(k_1,m') \in C_{1,+}(n)$. Repeating this procedure, we arrive at $\widehat u(n) = 0$. \qed

\bigskip
 
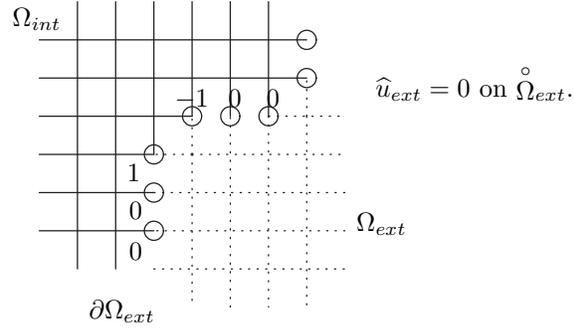
\begin{figure}[t]
\unitlength 0.1in
\begin{picture}( 19.0000, 16.0000)(  0.6000,-18.0000)
%
{\color[named]{Black}{%
\special{pn 8}%
\special{pa 1600 600}%
\special{pa 400 600}%
\special{fp}%
\special{pa 400 400}%
\special{pa 1600 400}%
\special{fp}%
\special{pa 1400 200}%
\special{pa 1400 800}%
\special{fp}%
\special{pa 1200 200}%
\special{pa 1200 800}%
\special{fp}%
\special{pa 1000 200}%
\special{pa 1000 800}%
\special{fp}%
\special{pa 800 200}%
\special{pa 800 800}%
\special{fp}%
\special{pa 600 200}%
\special{pa 600 1600}%
\special{fp}%
\special{pa 400 800}%
\special{pa 800 800}%
\special{fp}%
\special{pa 800 1000}%
\special{pa 400 1000}%
\special{fp}%
\special{pa 400 1200}%
\special{pa 800 1200}%
\special{fp}%
\special{pa 400 1400}%
\special{pa 800 1400}%
\special{fp}%
\special{pa 400 400}%
\special{pa 200 400}%
\special{fp}%
\special{pa 200 600}%
\special{pa 400 600}%
\special{fp}%
\special{pa 400 800}%
\special{pa 200 800}%
\special{fp}%
\special{pa 200 1000}%
\special{pa 400 1000}%
\special{fp}%
\special{pa 400 1200}%
\special{pa 200 1200}%
\special{fp}%
\special{pa 200 1400}%
\special{pa 400 1400}%
\special{fp}%
\special{pa 400 1600}%
\special{pa 400 200}%
\special{fp}%
}}%
%
{\color[named]{Black}{%
\special{pn 8}%
\special{pa 1000 1000}%
\special{pa 800 1000}%
\special{dt 0.045}%
\special{pa 1000 1000}%
\special{pa 1000 800}%
\special{dt 0.045}%
\special{pa 1000 1000}%
\special{pa 1800 1000}%
\special{dt 0.045}%
\special{pa 800 1200}%
\special{pa 1800 1200}%
\special{dt 0.045}%
\special{pa 1800 1400}%
\special{pa 800 1400}%
\special{dt 0.045}%
\special{pa 1000 1000}%
\special{pa 1000 1800}%
\special{dt 0.045}%
\special{pa 1200 1800}%
\special{pa 1200 800}%
\special{dt 0.045}%
\special{pa 1400 800}%
\special{pa 1400 1800}%
\special{dt 0.045}%
\special{pa 800 1600}%
\special{pa 1800 1600}%
\special{dt 0.045}%
\special{pa 1600 800}%
\special{pa 1600 1800}%
\special{dt 0.045}%
\special{pa 1600 800}%
\special{pa 1600 600}%
\special{dt 0.045}%
\special{pa 1400 800}%
\special{pa 1800 800}%
\special{dt 0.045}%
}}%
%
{\color[named]{Black}{%
\special{pn 8}%
\special{ar 1600 600 50 50  0.0000000 6.2831853}%
}}%
%
{\color[named]{Black}{%
\special{pn 8}%
\special{ar 1600 400 50 50  0.0000000 6.2831853}%
}}%
%
{\color[named]{Black}{%
\special{pn 8}%
\special{ar 1400 800 50 50  0.0000000 6.2831853}%
}}%
%
{\color[named]{Black}{%
\special{pn 8}%
\special{ar 1200 800 50 50  0.0000000 6.2831853}%
}}%
%
{\color[named]{Black}{%
\special{pn 8}%
\special{ar 1000 800 50 50  0.0000000 6.2831853}%
}}%
%
{\color[named]{Black}{%
\special{pn 8}%
\special{ar 800 1000 50 50  0.0000000 6.2831853}%
}}%
%
{\color[named]{Black}{%
\special{pn 8}%
\special{ar 800 1200 50 50  0.0000000 6.2831853}%
}}%
%
{\color[named]{Black}{%
\special{pn 8}%
\special{ar 800 1400 50 50  0.0000000 6.2831853}%
}}%
\put(6.6000,-11.4000){\makebox(0,0)[lb]{$1$}}%
\put(9.1000,-7.6000){\makebox(0,0)[lb]{$-1$}}%
\put(13.9000,-7.5000){\makebox(0,0)[lb]{$0$}}%
\put(6.7000,-13.4000){\makebox(0,0)[lb]{$0$}}%
\put(6.7000,-15.5000){\makebox(0,0)[lb]{$0$}}%
\put(11.9000,-7.5000){\makebox(0,0)[lb]{$0$}}%
\put(0.6000,-3.5000){\makebox(0,0)[lb]{$\Omega_{int}$}}%
\put(4.6000,-18.8000){\makebox(0,0)[lb]{$\partial \Omega_{ext}$}}%
\put(18.6000,-14.2000){\makebox(0,0)[lb]{$\Omega_{ext}$}}%
\put(19.6000,-7.2000){\makebox(0,0)[lb]{$\widehat{u}_{ext} =0$ on $\stackrel{\circ}\Omega_{ext}$.}}%
%
{\color[named]{Black}{%
\special{pn 8}%
\special{pa 800 800}%
\special{pa 800 1000}%
\special{fp}%
\special{pa 800 800}%
\special{pa 1000 800}%
\special{fp}%
}}%
\end{picture}%

\caption{A counter example for unique continuation property.}
\label{Fig2}
\end{figure}

An example of the domain which does not satisfy the cone condition is the one (in 2-dimension) whose boundary in the 4th quadrant has the form illustrated in Figure \ref{Fig2}, and is rectangular in the other quadrants. In this case, $\widehat u$ defined as in the figure satisfies 
$$
\big(\widehat H_{ext} -\frac{1}{2}\big)\widehat u = 0, \quad {\rm in} \quad \stackrel{\circ}\Omega_{ext},
$$
and $\widehat u = 0$ on $\stackrel{\circ}\Omega_{ext}$, however $\widehat u \not\equiv 0$ on $\partial\Omega_{ext}$.

\medskip


\subsection{Exterior eigenvalue problem}
Now let $\widehat H_{ext}^{(D)}$ be $\widehat H_{ext}$  subject to the Dirichlet boundary condition
\begin{equation}
\widehat{u} (n)=0 \quad \text{for any} \quad n\in \partial \Omega_{ext}, 
\label{S1DiricletOp}
\end{equation}
and $\widehat H_{ext}^{(R)}$  the Robin boundary condition
\begin{equation}
(\partial_{\nu} \widehat{u} )(n) +c(n) \widehat{u} (n) =0 \quad \text{for any} \quad n\in \partial \Omega_{ext},
\label{S1RobinOp}
\end{equation}
where $c(n)$ is a bounded real function on $\partial \Omega_{ext} $.
Here we deal with $\widehat{H}_{ext}^{(D)} $ as the bounded self-adjoint operator on $\ell^2 (\stackrel{\circ}\Omega_{ext} )$ with the Dirichlet boundary condition (\ref{S1DiricletOp}).
On the other hand, $\widehat{H}_{ext}^{(R)} $ is the bounded self-adjoint operator on $\ell^2 (\Omega_{ext} )$, interpreting $-\Delta_{disc} +\widehat{V} $ with (\ref{S1RobinOp}) as the Laplacian on graphs in view of (\ref{graph_laplacian}).

\begin{lemma}
$\sigma_{ess}(\widehat{H}_{ext}^{(D)}) = \sigma_{ess}( \widehat{H}_{ext}^{(R)}) = [0,d].$
\label{s2_lem_ess_spec}
\end{lemma}

Since this follows from the standard perturbation theory, we omit the proof.

\medskip

Theorem \ref{UniqueConti} asserts the non-existence of embedded eigenvalues for these operators. In particular, it is clear from the proof that for the Dirichlet case, because we have only to assume the cone condition for $\stackrel{\circ}\Omega_{ext}$.


\begin{theorem} \label{NonExistenceExterior}
(1) Let $\widehat H_{ext}^{(R)}$ be a Schr{\"o}dinger operator in an exterior domain $\Omega_{ext}$ with compactly supported potential subordinate to the Robin boundary condition.  Then $\sigma_p(\widehat H_{ext}^{(R)})\cap (0,d) = \emptyset$, if $\Omega_{ext}$ satisfies the cone condition. \\
\noindent
(2) Let $\widehat H_{ext}^{(D)}$ be a Schr{\"o}dinger operator in an exterior domain $\Omega_{ext}$ with compactly supported potential subordinate to the Dirichlet boundary condition.  Then $\sigma_p(\widehat H_{ext}^{(D)})\cap (0,d) = \emptyset$, if for any $n \in \; \stackrel{\circ}\Omega_{ext}$, there is a cone $C_{i,+}(n)$ (or $C_{i,-}(n)$) such that  $C_{i,+}(n) \subset \Omega_{ext}$ (or $C_{i,-}(n) \subset \Omega_{ext}$).
\end{theorem}


\section{Sobolev and Besov spaces on compact manifolds}

The condition (\ref{S1DecayCond}) is reformulated as a spectral condition for the Laplacian on the torus, which can further be rewritten by the Fourier transform. We do it on a compact Riemannian manifold in this section. 


\subsection{General case}
 
Let $M $ be a compact Riemannian manifold of dimension $d$ with the Riemannian metric $g $ and $L$ be the Laplace-Beltrami operator on $M$ defined by 
$$
L = - \sum_{i,j=1}^{d}\dfrac{1}{\sqrt{g}}\frac{\partial}{\partial x_i}\Big(\sqrt{g}g^{ij}\frac{\partial}{\partial x_j}\Big) , 
$$ 
where $ \big( g^{ij} \big)_{i,j=1}^d = g^{-1}$ and $\sqrt{g} = \sqrt{ \det g }$.

We introduce the Sobolev and Besov spaces on $M$ by two different ways.
One way is to use functions of $L$.
For $s \in {\bf R}$, we define $\mathcal H^s $ to be the completion of $C^{\infty}(M)$ by the norm
$\|\langle L\rangle^{s/2}u\|_{L^2 (M)}$.
We also define $\mathcal B^{\ast}$ to be the completion of $C^{\infty}(M)$ by the norm
$$
\left( \sup_{R>1} \dfrac{1}{R} \big\| \chi_R(\sqrt{L})u \big\|^2_{L^2 (M)} \right)^{1/2} ,
$$
where $\chi_R(t) = 1$ for $t < R$, $\chi_R(t) = 0$ for $t > R$.

Another way is to use the Fourier transform.
Let $\{\chi_j\}_{j=1}^N$ be a partition of unity on $M$ such that on each support of $\chi_j$, we can take one coordinate patch. 
We  define $H^s$ to be the completion of $C^{\infty}(M)$ by the norm
$ \big( \sum_{j=1}^N\|\langle\xi\rangle^s(\mathcal F\chi_ju)(\xi)\|^2_{L^2 ({\bf R}^d )} \big)^{1/2} $,
where $\mathcal F v = \widetilde v$ denotes the Fourier transform of $v$.
We define $B^{\ast}$ to be the completion of $C^{\infty}(M)$ by the norm
$$
\left( \sum_{j=1}^N\sup_{R>1}\frac{1}{ R} \big\| \chi_R( |\xi | ) ( \mathcal F \chi_j u) (\xi ) \big\| _{L^2 ({\bf R}^d )} ^2 \right) ^{1/2} .$$ 
The following inclusion relations hold for $s > 1/2$:
\begin{equation}
L^2 \subset \mathcal H^{-1/2} \subset \mathcal B^{\ast} \subset \mathcal H^{-s},
\quad 
L^2 \subset  H^{-1/2} \subset B^{\ast} \subset H^{-s}.
\label{AppendInclusion}
\end{equation}

These definitions of Sobolev and Besov spaces coincide. We show

\begin{lemma}\label{Hs=HsBs=Bs}
$\mathcal H^s = H^s$ for any $s \in {\bf R}$, and $\mathcal B^{\ast} = B^{\ast}.$
\end{lemma}

Proof. 
We prove $\mathcal B^{\ast} = B^{\ast}$.  It is well-known that $\mathcal H^s = H^s$ for $s \in {\bf R}$, whose proof is similar to, actually easier than, that for $\mathcal B^{\ast} = B^{\ast}$ given below. 

First let us recall a formula from functional calculus. Let $\psi(x) \in C^{\infty}({\bf R})$ be such that
\begin{equation}
|\psi^{(k)}(x)| \leq C_k\langle x\rangle^{m-k}, \quad \forall k \geq0,
\label{psi(x)estimates}
\end{equation}
for some $m \in {\bf R}$. One can construct  $\Psi(z) \in C^{\infty}({\bf C})$, called an {\it almost analytic extension} of $\psi$, having the following properties:
\begin{equation}
\left\{
\begin{split}
& \Psi(x) = \psi(x), \quad \forall x \in {\bf R}, \\
& |\Psi(z)| \leq C\langle z\rangle^m, \quad \forall z \in {\bf C}, \\
& |\overline{\partial_z}\Psi(z)| \leq C_n|{\rm Im}\,z|^{n}
\langle z\rangle^{m-n-1}, \quad \forall n \geq 1, \quad \forall z \in {\bf C},\\
& {\rm supp}\,\Psi(z) \subset \{z \, ; \, |{\rm Im}\,z| \leq 2 + 2|{\rm Re}\,z|\}.
\end{split}
\right.
\label{EstimateAlmostanalytic}
\end{equation}
In particular, if $\psi(x) \in C_0^{\infty}({\bf R})$, one can take $\Psi(z) \in C_0^{\infty}({\bf C})$.
Then, if $m < 0$,  for any self-adjoint operator $A$,  we have the following  formula
\begin{equation}
\psi(A) = \frac{1}{2\pi i}\int_{{\bf C}}\overline{\partial_z}\Psi(z)(z - A)^{-1}dzd\overline{z},
\label{FomulaHelfferSjostrand}
\end{equation}  
which is called the {\it formula of Helffer-Sj{\"o}strand}.
See \cite{HeSj}, \cite{DeGe}, p. 390.

We use a  semi-classical analysis employing $\hbar = 1/R$ as a small parameter (see e.g. \cite{Ro87}). We show that $\psi(\hbar^2 L)$ is equal to, modulo a lower order term, a pseudo-differential operator ($\Psi DO$) with symbol $\psi(\ell(x,\hbar \xi))$, where $\ell(x,\xi) = \sum_{i,j=1}^d g^{ij}(x)\xi_i\xi_j$. In fact, take $\chi(x), \chi_0(x) \in C^{\infty}(M)$ with small support such that $\chi_0(x) = 1$ on ${\rm supp}\,\chi$.  Consider a $\Psi DO$ $P_{\hbar}(z)$ with symbol $(\ell(x,\hbar \xi)-z)^{-1}$. Then 
$$
(\hbar^2L - z)\chi_0 P_{\hbar }(z)\chi = \chi +  Q_{\hbar}(z)\chi,
$$
where $Q_{\hbar}(z)$ is a $\Psi DO$ with symbol
\begin{equation}
\sum_{i=1}^2\hbar^i
\sum_{j=1}^3\frac{q_{ij}(x,\hbar\xi)}{(\ell(x,\hbar\xi) - z)^{j}}, \quad q_{ij}(x,\xi) \in S^{2j-1},
\label{SymbolQ(z)}
\end{equation}
$S^{m}$ being the standard H{\"o}rmander class of symbols (see \cite{HoVol3}, p. 65). This implies
$$
(\hbar^2L - z)^{-1}\chi = \chi_0 P_{\hbar}(z)\chi - (\hbar^2L-z)^{-1}Q_{\hbar}(z)\chi.
$$
By the symbolic calculus, we have
 \begin{equation}
\| \langle \hbar^2L \rangle ^{(s+1)/2} Q_{\hbar}(z)\chi \langle \hbar ^2 L \rangle ^{-s/2} \|_{{\bf B} (L^2 (M))} \leq\hbar \, C_{s,d}\,|{\rm Im}\, z|^{-N}\langle z\rangle^N,
\label{qz} 
\end{equation}
where $N > 0$ is a constant depending on $s$ and $d$, and $C_{s,d}$ does not depend on $\hbar$.
 
We take $\psi \in C^{\infty}_0({\bf R})$  and apply (\ref{FomulaHelfferSjostrand}). Then we have
\begin{equation}
\psi(\hbar^2L)\chi = \chi_0 \Psi_{\hbar}\chi + \Psi_{Q,\hbar}\chi,
\label{varphi(L)=}
\end{equation}
where $\Psi_{\hbar}$ is a $\Psi DO$ with symbol $\psi(\ell(x,\hbar\xi))$ and
\begin{equation}
\Psi_{Q,\hbar} = \frac{1}{2\pi i}\int_{\bf C}\overline{\partial_z}\Psi(z)(\hbar^2L-z)^{-1}Q_{\hbar}(z)dzd\overline{z}.
\label{PhiQ}
\end{equation}
We then have, letting $z = x+iy$,
\begin{gather*} 
\begin{split}
& \| \langle \hbar^2L \rangle ^{(s+3)/2} \Psi _{Q,\hbar} \chi\langle\hbar^2 L \rangle ^{-s/2} \| _{{\bf B} (L^2 (M))} \\
&\leq C\int_{ {\bf C} } | \overline{\partial _{z} } \Psi (z)| \| \langle\hbar^2 L \rangle ^{(s+3)/2} (\hbar^2L-z)^{-1} Q_{\hbar}(z)\chi \langle \hbar^2L \rangle ^{-s/2 } \| _{{\bf B} (L^2 (M))} dzd\overline{z} \\
&\leq C\int_{{\bf C}} |\overline{\partial _z} \Psi (z)| \| \langle \hbar^2 L \rangle ^{(s+3)/2}  (\hbar^2 L-z)^{-1} \langle \hbar^2 L \rangle ^{-(s+1)/2} \| _{{\bf B} (L^2 (M))}\\
& \ \ \ \ \ \ \ \ \ \ \ \ \ \ \  \ \ \ \ \ \ \ \ \ \  \ \ \ \ \ \ \ \ \cdot \| \langle \hbar^2 L \rangle ^{(s+1)/2}Q_{\hbar}(z)\chi\langle \hbar^2L\rangle^{-s/2} \|_{{\bf B} (L^2 (M))} dzd\overline{z}. 
\end{split} 
\end{gather*}
Since 
\begin{equation*}
\| \langle \hbar^2 L \rangle ^{(s+3)/2}  (\hbar^2 L-z)^{-1} \langle \hbar^2 L \rangle ^{-(s+1)/2}  \| _{{\bf B} (L^2 (M))} \leq \sup_{\lambda \in {\bf R}} \frac{\langle\lambda\rangle}{|\lambda -z |} \leq |\mathrm{Im}  \, z| ^{-1}  \langle z \rangle, \end{equation*}
we obtain, using (\ref{qz}),
\begin{equation} 
\begin{split} 
 \| \langle \hbar^2L & \rangle^{(s+3)/2} \Psi _{Q,\hbar} \chi\langle \hbar^2L \rangle ^{-s/2}\| _{{\bf B} (L^2 (M))} \\
& \leq \hbar C_{s,d} \int_{{\bf C}} | \overline{\partial _{z} } \Psi (z)| \ |\mathrm{Im} \, z |^{-1-N}\langle z\rangle^{N+1} dzd\overline{z} \\
&\leq \hbar C_{s,d}   \int_{{\bf C}} \langle z \rangle ^{m -1 } dzd\overline{z} , 
\end{split} 
\label{EstimatePsiQhbar}
\end{equation}
where we have used (\ref{EstimateAlmostanalytic}) with $m < -1 , n = N+1$.  This estimate implies
\begin{equation}
\sup_{\hbar < 1}\hbar \|\Psi_{Q,\hbar}\chi u\| ^2 _{L^2 (M)} < \infty, \quad 
{\rm if} \quad u \in \bigcap_{s>1/2} H^{-s}.
\label{PsiQhbaruOK}
\end{equation}
In fact, taking $0 \leq t \leq 3/2$, we have choosing $s = -3$ in (\ref{EstimatePsiQhbar}),
$$
\hbar^{1/2}\|\Psi_{Q,\hbar}\chi u\| _{L^2 (M)} \leq   C\hbar^{3/2}\|\langle \hbar^2L\rangle^{-t}u\| _{L^2 (M)} \leq
C\hbar^{-2t+3/2}\|\langle L\rangle^{-t}u\| _{L^2 (M)} .
$$
The right-hand side is bounded if $1/4 < t \leq 3/4$.

Now, by the definition of ${\mathcal B}^{\ast}$, we have the following equivalence 
\begin{equation}
u \in \mathcal B^{\ast} \Longleftrightarrow
\left\{
\begin{split}
& u \in H^{-s}, \quad \forall s > 1/2, \\
& \sup_{\hbar<1}\hbar \big\| \psi(\hbar^2L)u \big\| ^2 _{L^2 (M)} < \infty, \quad \forall \psi \in C_0^{\infty}({\bf R}). 
\end{split}
\right.
\label{uinBastrewrite}
\end{equation}
In fact, the left-hand side is equivalent to the right-hand side for one fixed $\psi$ such that $\psi(t) = 1$ for $|t|<1$, and $\psi(t) = 0$ for $|t|>2$.

By virtue of (\ref{varphi(L)=}) and (\ref{PsiQhbaruOK}), (\ref{uinBastrewrite}) is equivalent to
\begin{equation}
\left\{
\begin{split}
& u \in H^{-s}, \quad \forall s > 1/2, \\
& \sup_{\hbar<1}\hbar \big\|\chi_0\Psi_{\hbar}\chi_j u \big\|^2 _{L^2 (M)} < \infty, \quad \forall j. 
\end{split}
\right.
\label{AppendAEquiv2}
\end{equation}
The symbol of $(\Psi_{\hbar})^{\ast}\chi_0(x)^2\Psi_{\hbar}$ is equal to 
\begin{equation*}
\chi_0 (x)^2 \psi ( \ell (x,\hbar \xi ) )^2 + O(\hbar ) .
\end{equation*}
Then by a suitable choice of $0 < c_1 < c_2$, we have
$$
\psi \Big(\frac{\hbar ^2 |\xi |^2}{c_1} \Big)^2  \leq \psi \left(\ell(x,\hbar\xi)\right) ^2 \leq 
\psi \Big( \frac{\hbar ^2 | \xi|^2 }{c_2 }\Big)^2.
$$
Moreover, we can assume that there exists $q(x,\xi) \in C^{\infty}( M \times{\bf R}^d)$ such that
\begin{equation}
\left\{
\begin{split}
& \psi (|\xi|^2 /c_2) ^2 - \psi(\ell(x,\xi))^2 = q(x,\xi)^2, \\
& {\rm supp}\ q(x,\xi) \subset M \times\{a < |\xi|< b\},
\end{split}
\right.
\label{Appndqxxisupport}
\end{equation}
for some $0 < a < b$. Since the symbol of $\psi\Big(\frac{-\hbar^2\Delta}{c_2 }\Big)\chi_0(x)^2\psi\Big(\frac{-\hbar^2\Delta}{c_2 }\Big) $ is equal to 
\begin{equation*}
\chi_0 (x)^2 \psi \Big( \frac{\hbar^2 |\xi |^2 }{c_2  } \Big)^2 +O(\hbar), 
\end{equation*}
we see that the symbol of $\psi\Big(\frac{-\hbar^2\Delta}{c_2 }\Big)\chi_0(x)^2\psi\Big(\frac{-\hbar^2\Delta}{c_2  }\Big) - 
\Psi_{\hbar}^{\ast}\chi_0(x)^2\Psi_{\hbar}$ is estimated as
\begin{gather*}
\begin{split}
& \chi_0 (x)^2 \psi \Big( \frac{\hbar^2|\xi|^2 }{c_2 } \Big) ^2 -\chi_0 (x)^2 \psi \Big( \ell (x,\hbar \xi ) \Big) ^2  +O(\hbar) \\
=& \chi_0 (x)^2 q(x,\hbar\xi )^2 +O(\hbar ).
\end{split}
\end{gather*}
Then we have
\begin{equation}
\begin{split}
& \psi\Big(\frac{-\hbar^2\Delta}{c_2  }\Big)\chi_0(x)^2\psi\Big(\frac{-\hbar^2\Delta}{c_2  }\Big) - 
(\Psi_{\hbar})^{\ast}\chi_0(x)^2\Psi_{\hbar} \\
 =&( Q_{0,\hbar} )^{\ast}\, \chi_0 (x)^2 Q_{0,\hbar} + Q_{1,\hbar}.
\end{split}
\label{Psihchi0Psihfromabove}
\end{equation}
In the right-hand side, $Q_{0,\hbar}$ is a $\Psi DO$ with symbol $q(x,\hbar\xi)$, where $q(x,\xi )$ is given by (\ref{Appndqxxisupport}), and $Q_{1,\hbar}$ is a $\Psi DO$ with symbol $q_1(x,\xi;\hbar)$ admitting the asymptotic expansion 
\begin{equation}
q_1(x,\xi;\hbar) \sim \sum_{j\geq 1}\hbar^{j} q_j (x,\hbar\xi  ),
\label{q1xxiasymptoticexpand}
\end{equation}
with $q_j (x,\xi )$ having the same support property as in (\ref{Appndqxxisupport}). Since the 1st term of the right-hand side of (\ref{Psihchi0Psihfromabove}) is non-negative, we have proven
\begin{equation}
\begin{split}
 \psi\Big(\frac{-\hbar^2\Delta}{c_2  }\Big)\chi_0(x)^2\psi\Big(\frac{-\hbar^2\Delta}{c_2  }\Big) \geq 
(\Psi_{\hbar})^{\ast}\chi_0(x)^2\Psi_{\hbar} + Q_{1,\hbar}.
\end{split}
\label{Fromabove}
\end{equation}
 By a similar computation, we can prove
\begin{equation}
\begin{split}
 (\Psi_{\hbar})^{\ast}\chi_0(x)^2\Psi_{\hbar} \geq \psi\Big(\frac{-\hbar^2\Delta}{c_1}\Big)\chi_0(x)^2\psi\Big(\frac{-\hbar^2\Delta}{c_1 }\Big) 
 + Q_{1,\hbar}'.
\end{split}
\label{Frombelow}
\end{equation}
 By (\ref{Fromabove}) and (\ref{Frombelow}), we have
\begin{equation}
\begin{split}
 \hbar\psi\Big(\frac{-\hbar^2\Delta}{c_2  }\Big)\chi_0(x)^2\psi\Big(\frac{-\hbar^2\Delta}{c_2  }\Big) -  \hbar Q_{1,\hbar} & \geq 
\hbar(\Psi_{\hbar})^{\ast}\chi_0(x)^2\Psi_{\hbar}  \\
& \geq \hbar\psi\Big(\frac{-\hbar^2\Delta}{c_1}\Big)\chi_0(x)^2\psi\Big(\frac{-\hbar^2\Delta}{c_1 }\Big) 
 + \hbar \widetilde Q_{1,\hbar}',
\end{split}
\label{Fromabovebelow}
\end{equation}
where $Q_{1,\hbar}$ and $Q_{1,\hbar}'$ have the property (\ref{q1xxiasymptoticexpand}). As $u \in H^{-s}, \ \forall s >1/2$, we have
\begin{equation}
\sup_{\hbar<1}\hbar\big|(Q_{1,\hbar}u,u)\big| + \sup_{\hbar<1}\hbar\big|(Q_{1,\hbar}'u,u)\big|< \infty.
\nonumber
\end{equation}
Therefore, by (\ref{Fromabovebelow}), 
$\displaystyle{
\sup_{\hbar<1}\hbar\|\chi_0\Psi_{\hbar}\chi_j u\|^2 < \infty}$
is equivalent to 
$$
\sup_{R>1}\frac{1}{R}\int_{{\bf R}^d} \Big|\psi \Big (\frac{|\xi |^2 }{c R^2 } \Big) \widehat{\chi _j u } (\xi ) \Big|^2 d\xi < \infty,
$$
for some $c > 0$,
which is equivalent to $u \in B^{\ast}$. We have thus completed the proof of Lemma \ref{Hs=HsBs=Bs}. \qed

\medskip
By the same argument, we can also prove the following lemma. 


\begin{lemma}
If $u \in \mathcal B^{\ast}$, we have the following equivalence  
$$
\lim_{R\to \infty} \frac{1}{\sqrt{R}}\Big\| \psi\Big(\frac{ L }{R^2 }\Big)u \Big\| _{L^2 (M)} = 0 \Longleftrightarrow 
\lim_{R\to \infty} \frac{1}{\sqrt R}\Big\| \psi\Big(\frac{|\xi|^2}{R^2} \Big)\big(\mathcal F\chi_j u \big)(\xi) \Big\| _{L^2 ({\bf R}^d ) } = 0 ,
$$
for any $j$ and any $\psi \in C_0^{\infty}({\bf R})$, where $\{\chi_j\}_{j=1}^N$ is the partition of unity on $M$.
\label{S3eq1}
\end{lemma}


\subsection{Torus}
We interpret the above results for the case of the torus ${\bf T}^d $ defined by (\ref{S1Torus}).
Let ${\mathcal U}$ be the unitary operator from
$\ell^{2}({\bf Z}^d)$ to $L^{2}({\bf T}^d)$ defined by
\begin{equation}
({\mathcal U}\,{\widehat f})(x) = (2\pi)^{-d/2}\sum_{n\in{\bf Z}^d}{\widehat f}(n)
e^{-in\cdot x}.
\label{S2fourier}
\end{equation}
Letting
\begin{equation}
H_0 = {\mathcal U}\, {\widehat H_0}\, {\mathcal U}^{\ast}, \quad \widehat{H}_0 = -\Delta_{disc},
\nonumber
\end{equation}
we have
\begin{equation}
H_0 = h(x)= \frac{1}{2}\Big(d - \sum_{j=1}^d\cos x_j\Big) =\sum_{j=1}^d \sin^2 \Big( \frac{x_j }{2} \Big) ,
\label{S2H0cosine}
\end{equation}
We define operators $\widehat N_j$ and $N_j$ by
$$
\big(\widehat N_j\widehat f)(n) = n_j\widehat f(n),
\quad
N_j = \mathcal U\widehat N_j\mathcal U^{\ast} = i\frac{\partial}{\partial x_j}.
$$
We put $N = (N_1,\cdots,N_d)$, and let $N^2$ be the self-adjont operator defined by
\begin{equation}
N^2 = \sum_{j=1}^dN_j^2 =  - \Delta, \quad {\rm on} \quad {\bf T}^d,
\nonumber
\end{equation}
where $\Delta$ denotes the Laplacian on ${\bf T}^d = [-\pi,\pi]^d$ with periodic boundary condition. We can then apply the results in the previous subsection 
to $L = - \Delta$. We put
\begin{equation}
|N| = \sqrt{N^2} = \sqrt{-\Delta}.
\nonumber
\end{equation}
In the following, we simply denote $ \| \, \cdot \, \|_{L^2 ({\bf T}^d )} = \| \, \cdot \, \| $. 
For $s \in {\bf R}$,
let ${\mathcal H}^s$ be the completion of $D(|N|^s)$ with respect to the norm
$\|u\|_{s} = \|\langle N\rangle^{s}u\|$ i.e.
\begin{equation}
{\mathcal H}^s = \big\{u \in {\mathcal D}^{\prime}({\bf T}^d)\, ; \, \|u\|_{s} = \|\langle N\rangle^{s}u\| <
\infty \big\},
\nonumber
\end{equation}
where $ \mathcal{D}' ({\bf T}^d ) $ denotes the space of distribution on ${\bf T}^d $.
Put $\mathcal H = \mathcal H^0 = L^2({\bf T}^d)$.

For a self-adjoint operator $T$, let $\chi(a \leq T < b)$ denote the operator $\chi_{I}(T)$, where $\chi_I(\lambda)$ is the characteristic function of the interval $I = [a,b)$. The operators $\chi(T < a)$ and $\chi(T \geq b)$ are defined similarly.  Using the series $\{r_j\}_{j=0}^{\infty}$ with $r_{-1} = 0$, $r_j = 2^j \ (j \geq 0)$, we define the Besov space $\mathcal B$ by
\begin{equation}
\mathcal B = \Big\{f \in {\mathcal H}\, ; \|f\|_{\mathcal B} = \sum_{j=0}^{\infty}r_j^{1/2}\|\chi(r_{j-1} \leq |N| < r_j)f\| < \infty\Big\}.
\nonumber
\end{equation}
Its dual space $\mathcal B^{\ast}$ is the completion of $\mathcal H$ by the following norm
\begin{equation}
\|u\|_{\mathcal B^{\ast}}= \sup_{j\geq 0}2^{-j/2}\|\chi(r_{j-1} \leq |N| < r_j)u\|.\nonumber
\end{equation}
The following Lemma \ref{S3AgHo} is proved in the same way as in  \cite{AgHo76}.


\begin{lemma}
(1) There exists a constant $C > 0$ such that
\begin{equation}
C^{-1}\|u\|_{\mathcal B^{\ast}}\leq \left(\sup_{R>1}\frac{1}{R}\|\chi(|N| < R)u\|^2\right)^{1/2}\leq C\|u\|_{\mathcal B^{\ast}}.
\nonumber
\end{equation}
(2) For $s > 1/2$, the following inclusion relations hold:
\begin{equation}
\mathcal H^{s} \subset \mathcal B \subset \mathcal H^{1/2} \subset \mathcal H \subset \mathcal H^{-1/2} \subset \mathcal B^{\ast} \subset \mathcal H^{-s}.
\nonumber
\end{equation}
\label{S3AgHo}
\end{lemma}

In view of the above lemma, in the following, we use
\begin{equation}
\|u\|_{\mathcal B^{\ast}} = \left(\sup_{R>1}\frac{1}{R}\|\chi(|N| < R)u\|^2\right)^{1/2}
\nonumber
\end{equation}
as a norm on $\mathcal B^{\ast}$.

We also put $\widehat{\mathcal H} = \ell^2({\bf Z}^d)$, and define $\widehat{\mathcal H}^s$, $\widehat{\mathcal B}$, $\widehat{\mathcal B}^{\ast}$ by replacing $N$ by $\widehat N$. Note that $\widehat{\mathcal H}^s = \mathcal U^{\ast}\mathcal H^s$ and so on. In particular, 
Parseval's formula implies that
\begin{equation}
\|u\|_{\mathcal H^s}^2 = \|\widehat u\|_{\widehat{\mathcal H}^s}^2 = \sum_{n\in {\bf Z}^d}(1 + |n|^2)^s|\widehat u(n)|^2 ,
\nonumber
\end{equation}
\begin{equation}
\|u\|_{{\mathcal B}^{\ast}}^2 =  
\|\widehat u\|_{\widehat{\mathcal B}^{\ast}}^2
= 
\sup_{R>1}\frac{1}{R}\sum_{|n|<R}|\widehat u(n)|^2,
\nonumber
\end{equation}
$\widehat u(n)$ being the Fourier coefficient of $u(x)$.


\section{Proof of Theorem 1.1}
\subsection{Some remarks for the Fermi surface}





The Fermi surface $M_{\lambda}$ is not a smooth submanifold of ${\bf T}^d $ in general.
Here we consider some properties of $M_{\lambda }$ as an analytic set (see e.g. \cite{Chi}, \cite{KuVa}).

Let $ {\bf T}^d _{{\bf C}} = {\bf C}^d /(2\pi {\bf Z} )^d $ be the complex torus and define 
\begin{equation}
M_{\lambda }^{{\bf C} } =\big\{ z \in {\bf T}^d_{{\bf C}} \ ; \ h(z)=\lambda \big\} . \label{isospectral_complex} 
\end{equation}
Then $M_{\lambda}^{{\bf C}} \cap {\bf R}^d = M_{\lambda}$.
  .

\begin{lemma}
(1) For $\lambda \in (0,d) \setminus {\bf Z} $, $M_{\lambda}^{{\bf C}}$ is a $(d-1)$-dimensional, closed submanifold of ${\bf T}^d_{{\bf C}}$.
\\
(2) For $\lambda \in (0,d) \cap {\bf Z} $, $M_{\lambda}^{{\bf C}}$ consists of a disjoint union $\big( \mathrm{reg}\, M_{\lambda}^{{\bf C}} \big) \cup \big( \mathrm{sng} \, M_{\lambda}^{{\bf C}} \big) $, where 
\begin{gather}
\mathrm{sng} \, M_{\lambda}^{{\bf C}} = \big\{ z\in M_{\lambda}^{{\bf C}} \ ; \ z_j \equiv 0 \ (\mathrm{mod} \ \pi ) \ \text{for all} \ j=1, \cdots ,d \big\} , \label{sing_part} \\
\mathrm{reg} \, M_{\lambda}^{{\bf C}} = M_{\lambda}^{{\bf C}} \setminus \big( \mathrm{sng} \, M_{\lambda}^{{\bf C}}  \big). \label{reg_part} 
\end{gather}
Moreover, $\mathrm{reg} \, M_{\lambda}^{{\bf C}} $ is a $(d-1)$-dimensional, open submanifold of ${\bf T}_{{\bf C}}^d $.
\label{s4_lem_analyticset}
\end{lemma}

Proof.
Since $\nabla h(z)= \frac{1}{2} ( \sin z_1 , \cdots , \sin z_d )$, $\nabla h(z)=0$ if and only if $z_j \equiv 0 \ (\mathrm{mod} \ \pi )$ for all $j=1, \cdots , d$.
For $\lambda \in (0,d)\setminus {\bf Z}$, the intersection $M_{\lambda}^{{\bf C}} \cap \{ z\in {\bf T}^d_{{\bf C}} \ ; \ z_j \equiv 0 \ (\mathrm{mod} \ \pi) \ \text{for all} \ j=1, \cdots ,d \} $ is empty, and for $\lambda \in (0,d)\cap {\bf Z}$, we see $ \mathrm{sng}\, M_{\lambda}^{{\bf C}} \not= \emptyset $.
By the definitions of $\mathrm{reg} \, M_{\lambda}^{{\bf C}} $ and $\mathrm{sng} \, M_{\lambda}^{{\bf C}} $, we see that $M_{\lambda}^{{\bf C}} = \big( \mathrm{reg}\, M_{\lambda}^{{\bf C}} \big) \cup \big( \mathrm{sng} \, M_{\lambda}^{{\bf C}} \big) $ as a disjoint union. \qed

\medskip

For $\lambda \in (0,d) \setminus {\bf Z}$, $M_{\lambda }^{{\bf C}} $ is irreducible i.e. a connected complex submanifold (see \cite{Sha}, proof of Theorem 9).
On the other hand, if $\lambda \in (0,d)\cap {\bf Z}$, $ \mathrm{reg}\, M_{\lambda}^{{\bf C}} $ may not be a connected submanifold.
In fact, for $d=2 $ and $\lambda =1$, we have $M_{\lambda }^{{\bf C}} =A_+ \cup A_-$ where
$$
A_{\pm} =\{ z\in {\bf T}_{{\bf C}}^2 \ ; \ z_2 \equiv \pm z_1 +\pi \ (\mathrm{mod} \ 2\pi ) \}.
$$
Moreover, $A_+ \cap A_- = \mathrm{sng} \, M_{\lambda}^{{\bf C}}$.

\medskip

Now we give slightly weaker statement than the irreducibility of $M_{\lambda}^{{\bf C}} $.
The following fact implies that any irreducible components of $\mathrm{reg} \, M_{\lambda}^{{\bf C}}$ intersects an open part of $M_{\lambda}$.

\begin{lemma}
Let $\lambda \in (0,d)$ and $ \{ S_{\lambda ,j} \}_j $ be the set of connected components of $\mathrm{reg}\, M_{\lambda}^{{\bf C}}$.
Then, for any $S_{\lambda ,j} $, the intersection $S_{\lambda , j} \cap M_{\lambda} $ contains an open part of $\mathrm{reg}\, M_{\lambda}^{{\bf C}} $.

\label{s4_lem_connect}
\end{lemma}

Proof. Using the change of variables ${\bf T}_{{\bf C}} \ni z_j \mapsto w_j \in \mathcal{X}^d$ which is the Riemann surface of $\arccos$, more precisely $ w_j = \arccos z_j $, $j=1,\cdots ,d $, we can reduce the equation $h(z)-\lambda =0$ on ${\bf T}^d_{{\bf C}} $ to $ w_1 +\cdots +w_d = d-2\lambda $.
Then we can see that any connected component of the algebraic variety 
$$
V_{\lambda} =\{ w\in \mathcal{X}^d \ ; \ w_1 +\cdots + w_d =d-2\lambda \} 
$$ 
have some intersection containing an open part of $ V_{\lambda} \cap {\bf R}^d $. \qed


\subsection{Helmholtz equation}

We extend $\widehat{u} (n) $ to be zero for $|n|\leq R_0$ and denote it by $\widehat{u} $ again.
Then we have 
\begin{equation}
(\widehat{H}_0 -\lambda )\widehat{u} =\widehat{f}, 
\label{helmholtz} 
\end{equation}
where $\widehat f$ is compactly supported. In fact, letting ${\widehat P}(k)$ be the projection onto the site $k$, it is written as $\widehat f = \sum_{|k|\leq R_0+1}c_k\widehat P(k)\widehat u$. 

Here we derive the proof for $\lambda \in (0,d)\cap {\bf Z}$.
When $\lambda \in (0,d)\setminus {\bf Z}$, 
we can follow the same argument as $ \mathrm{sng}\, M_{\lambda}^{{\bf C}} =\emptyset $ and the proof is slightly easier than that for $\lambda \in (0,d)\cap {\bf Z}$.
We first note the following Lemma.


\begin{lemma}
Let $\lambda \in (0,d)$ and $\widehat{u} $ satisfy (\ref{S1Equation}) and (\ref{S1DecayCond}).
Then $u \in C^{\infty}({\bf T}^d \setminus (\mathrm{sng}\, M_{\lambda}^{{\bf C}}))$ and $\widehat{f} $  satisfies 
\begin{equation}
(\mathcal{U} \widehat{f} )(x)= f(x)=0 \quad {\rm  on} \quad (\mathrm{reg} \,M_{\lambda } ^{{\bf C}} )\cap {\bf T}^d .
\label{rellich_perturbation}
\end{equation} 

\label{lem_s4_fzero_real}
\end{lemma}

Proof. Passing to the Fourier series, (\ref{S1DecayCond}) implies 
\begin{equation}
\lim_{R \to \infty} \frac{1}{R} \int_{{\bf T}^d} | \chi (|N|<R ) u(x) |^2 dx=0 . \label{rellich_asymptoticzero} 
\end{equation}

Take a point $ x^{(0)}\in (\mathrm{reg}\, M_{\lambda}^{{\bf C}}) \cap {\bf T}^d $ and fix it.
Let $U$ be a sufficiently small neighborhood of $x^{(0)} $ in ${\bf T}^d $ such that $ U\cap (\mathrm{sng} \, M_{\lambda}^{{\bf C}} ) = \emptyset $.
Making the change of variables $x \mapsto (y_1 , y')$ so that $ y_1 =h(x)-\lambda $ and $y'=(y_2 ,\cdots ,y_d )$ in $U$, 
the Laplacian $N^2 =-\Delta $ on ${\bf T}^d $ is translated as the Laplace-Beltrami operator in $y$-coordinate using the Jacobian.
Letting $\chi \in C^{\infty } ({\bf T}^d )$ be such that $\chi (x^{(0)} )=1 $, $\mathrm{supp}\, \chi \subset U$, 
we have that, by Lemma \ref{S3eq1}, (\ref{rellich_asymptoticzero}) leads to
\begin{equation}
\lim_{R \to \infty} \frac{1}{R} \int_{|\eta |<R} |\widetilde{ \chi u} (\eta )|^2 d\eta =0, 
\label{rellich_asymptoticzero2} 
\end{equation} 
where $\widetilde{\chi u} (\eta )$ is the Fourier transform of $\chi u$:
$$
\widetilde{\chi u} (\eta )= (2\pi )^{-d/2} \int_{{\bf R}^d } e^{-iy\cdot \eta } \chi (y) u(y) dy .
$$

By (\ref{helmholtz}), $u$ satisfies
 \begin{equation}
(h(x)-\lambda )u=f \quad \text{on} \quad {\bf T}^d , \label{helmholtz2}
 \end{equation}
where  $f$ is a polynomial of $e^{ix_j } $, $j=1 , \cdots , d $, since $\widehat f$ is compactly supported.
Letting $ u_{\chi} (x)=\chi(x)u(x) $, $f_{\chi} (x)=\chi (x)f(x)$ and making the change of variable $x \mapsto y$ as above, we have by passing to the Fourier transform, $\frac{\partial}{\partial \eta _1} \widetilde{u_{\chi}} (\eta)= -i \widetilde{f_{\chi}} (\eta )$. 
Integrating this equation, we have 
$$
\widetilde{u_{\chi} }(\eta) = - i\int_0^{\eta_1}\widetilde{f_{\chi}}(s, \eta' )ds + \widetilde{u_{\chi}} (0, \eta' ), \quad \eta'=(\eta_2 ,\cdots , \eta_d ).
$$
Since $\widetilde{f_{\chi}} (\eta)$ is rapidly decreasing, we then see the existence of the limit
$$
\lim_{\eta_1\to\infty}\widetilde{u_{\chi}}(\eta) = - i \int_0^{\infty}\widetilde{f_{\chi}} (s, \eta' )ds + \widetilde{u_{\chi}} (0, \eta' ).
$$
We show that this limit vanishes. 
Let $D_R $ be the slab such that 
\begin{equation*}
D_R = \left\{ \eta \ ; \ |\eta '|<\delta R , \ \frac{R}{3 } < \eta_1 < \frac{2R}{3} \right\}. 
\end{equation*}
Then we have $D_R \subset \{ |\eta |<R \} $ for a sufficiently small $\delta >0 $.
We then see that 
\begin{equation*}
\frac{1}{R} \int_{ D_R } | \widetilde{u_{\chi}} (\eta )|^2 d\eta =\frac{1}{R} \int_{|\eta '|<\delta R} \int_{R/3}^{2R/3} |\widetilde{u_{\chi}} (\eta_1,\eta' )|^2 d\eta_1 d\eta'  \leq \frac{1}{R} \int_{ |\eta |<R} |\widetilde{u_{\chi}} ( \eta )|^2 d\eta. 
\end{equation*}
As $R\rightarrow \infty $, the right-hand side tends to zero by (\ref{rellich_asymptoticzero2}), hence so does the left-hand side, which proves that 
$\lim_{\eta_1\to\infty}\widetilde{u_{\chi} } (\eta) = 0$.

We have, therefore, 
\begin{equation}
\widetilde{u_{\chi}} (\eta )= i\int_{\eta_1 } ^{\infty } \widetilde{f_{\chi}} (s ,\eta ')ds.
\nonumber
\end{equation}
This shows that $u_{\chi} =\chi u \in C^{\infty } ({\bf T}^d \setminus ( \mathrm{sng}\, M_{\lambda}^{{\bf C}}))$.

It is easy to see that $u$ is smooth outside $M_{\lambda }$.
Then $u\in C^{\infty} ({\bf T}^d \setminus (\mathrm{sng}\, M_{\lambda}^{{\bf C}}))$.
In particular, $f (x) =0$ in $ (\mathrm{reg}\, M_{\lambda}^{{\bf C}} ) \cap {\bf T}^d $ by (\ref{helmholtz2}).
\qed

\medskip

In the following discussion, we use the function theory of several complex variables.
We extend $f(x)$ to the polynomial of $e^{iz_j}$ for $z_j \in {\bf C}$, $j=1, \cdots ,d$.
Lemmas \ref{s4_lem_connect} and \ref{lem_s4_fzero_real} imply that the zeros of $f$ can be extended to $ M_{\lambda}^{{\bf C}}$.
For the proof, see Corollary 7 of \cite{KuVa}.

\begin{lemma}
Let $\lambda \in (0,d)$.
We have $f(z)=0$ on $\mathrm{reg} \, M_{\lambda}^{{\bf C}}$.
\label{rellich_continuation} 
\end{lemma}

\medskip

Let us return to the equation (\ref{helmholtz}).
Take any domain $D \subset {\bf T}^d_{{\bf C}} \setminus (\mathrm{sng}\, M_{\lambda}^{{\bf C}} )$.
By Lemma \ref{rellich_continuation}, and using the Taylor expansion, we see that there exists a holomorphic function $g$ in the domain $D$ such that $f(z)= \big( h(z)-\lambda \big) g(z) $, hence $u(z)= f(z)/ \big( h(z)-\lambda \big)$ is an holomorphic function of $z \in D$, since $ (\partial _{z_1} h(z), \cdots , \partial _{z_d} h(z) ) \not= 0 $ in $D$.
The uniqueness theorem for holomorphic functions leads that $  f(z)/(h(z)-\lambda ) $ is holomorphic in $ {\bf T}^d_{{\bf C}} \setminus ( \mathrm{sng}\, M_{\lambda}^{{\bf C}} )$.
However, since we have assumed $d\geq 2$, and $\mathrm{sng}\, M_{\lambda}^{{\bf C}} $ is a $0$-dimensinal set, $\mathrm{sng}\, M_{\lambda}^{{\bf C}} $ is a set of removable singularities (from Hartogs' extension theorem or see e.g. Corollary 7.3.2 of \cite{Kr82}).
Then $  f(z)/(h(z)-\lambda )$ can be extended to an entire function on ${\bf T}_{{\bf C}} ^d$ uniquely.
We denote it $f(z)/(h(z)-\lambda )$ again.

Here we pass to the variables $w_j =e^{iz_j }$, $j=1, \cdots , d $.
Note that the map 
\begin{equation*}
{\bf T}_{{\bf C}} ^d \ni z \mapsto w \in {\bf C}^d \setminus \bigcup_{j=1}^d A_j , \quad A_j =\{ w \in {\bf C}^d \ ; \ w_j =0 \} , 
\end{equation*}
is biholomorphic.
Then there exist positive integers $\alpha _j $ such that 
\begin{equation*}
f(z)=\sum c_{\beta}w^{\beta}=F(w) \prod_{j=1}^d w_j ^{-\alpha_j }, 
\end{equation*}
where $F(w) $ is a polynomial, $F(w) = \sum a_{\gamma}w^{\gamma}$, with the property that
$$
a_{\gamma} = 0, \ {\rm if} \ {\rm one} \ {\rm of}\ \gamma_j \leq 0 \ {\rm in} \ 
\gamma = (\gamma_1,\cdots,\gamma_d). 
$$
We  factorize $h(z)-\lambda $ as 
\begin{equation*}
h(z)-\lambda = \frac{d}{2}-\lambda -\frac{1}{4} \sum_{j=1}^d (w_j + w_j^{-1} ) = H_{\lambda } (w) \prod_{j=1}^d w_j^{-1} , 
\end{equation*}
where 
\begin{equation}
H_{\lambda } (w)= \Big( \frac{d}{2} -\lambda \Big) \prod_{j=1}^d w_j -\frac{1}{4} \Big( \sum_{ j=1}^d w_j \Big) \prod_{j=1} ^d w_j - \frac{1}{4}\sum_{j=1}^d \Big( \prod_{i \not= j } w_i \Big). 
\label{AppBHlambdaw}
\end{equation}
Then 
\begin{equation}
\frac{f(z)}{h(z)-\lambda } =\frac{ F(w)}{H_{\lambda } (w) } \prod_{j=1}^d w_j^{1-\alpha_j } . \label{variableschange} 
\end{equation}
Since $f(z)/(h(z)-\lambda ) $ is analytic, $F(w)/H_{\lambda } (w) $ is also analytic  except possibly on hyperplanes $A_j $, $j=1, \cdots , d$.
However, due to the expression (\ref{AppBHlambdaw}), we have
\begin{equation}
H_{\lambda}(w) \neq 0, \ {\rm if} \ w \in \bigcup_{k=1}^d V_k,
\nonumber
\end{equation}
$$
V_k = \{(w_1,\cdots,w_{k-1},0,w_{k+1},\cdots,w_d)\, ; \, w_i \neq 0, i \neq k\}.$$
Hence $F(w)/H_{\lambda } (w) $ is analytic except only on some sets of complex dimension $d-2$ (the intersection of two hyperplanes).
Therefore, 
\begin{itemize}
\item $F(w)/H_{\lambda } (w) $ is an entire function. 
\end{itemize}
See e.g. Corollary 7.3.2 of \cite{Kr82} again.
In particular, 
\begin{itemize}
\item
 $ F(w)=0 $ on the set $\{ w\in {\bf C}^d \ ; \ H_{\lambda }(w)=0\} $.
\end{itemize}

Finally, we use the following fact, a corollary of the Hilbert Nullstellensatz (See e.g. Appendix 6 of \cite{Shafa}). Let ${\bf C}[w_1,\cdots,w_d]$ be the ring of polynomials of variables $w_1,\cdots,w_d$.


\begin{lemma}
If $f, g \in {\bf C}[w_1,\cdots,w_d] $, and suppose that $f$ is irreducible. If $g=0$ on all zeros of $f$, there exists $h \in {\bf C}[w_1,\cdots,w_d]$ such that $g=fh$.
\label{hilbertnullstellensatz} 
\end{lemma}

\medskip

We factorize $H_{\lambda}(w)$ so that
$$
H_{\lambda } (w)= H^{(1)} _{\lambda } (w) \cdots H^{(N)} _{\lambda } (w),
$$
where each $H_{\lambda}^{(j)} (w) $ is an irreducible polynomial.
We prove inductively that 
$$
F(w)/H_{\lambda }^{(1)}(w)\cdots H_{\lambda}^{(k)}(w)\ {\rm is}\ {\rm a}\ {\rm polynomial}\ {\rm for}\  1 \leq k \leq N.
$$
Note that, since we know already that $F(w)/H_{\lambda}(w)$ is entire, 
\begin{itemize}
\item $F(w)/H_{\lambda }^{(1)}(w)\cdots H_{\lambda}^{(k)}(w)$ is also entire, 
\item $F(w) = 0$ on the zeros of $H_{\lambda }^{(1)}(w)\cdots H_{\lambda}^{(k)}(w)$.
\end{itemize}

Consider the case $k=1$. Since  $F(w) = 0$ on the zeros of $H_{\lambda}^{(1)}(w)$,
 Lemma \ref{hilbertnullstellensatz} implies that  $F(w) / H_{\lambda}^{(1)} (w)$ is a polynomial.
 
Assuming the case $k \leq \ell -1$, we consider the case $k=\ell$.
By the induction hypothesis, there exists a polynomial $P_{\ell -1}(w)$ such that 
\begin{equation*} 
\frac{F(w)}{  H_{\lambda}^{(1)} (w) \cdots H_{\lambda}^{(\ell-1 )} (w) }=P_{\ell-1}(w) . 
\end{equation*}
Then we have $F(w)/ ( H_{\lambda}^{(1) } (w) \cdots H_{\lambda}^{(\ell)} (w) ) =P_{\ell-1}(w)/H^{(\ell)}_{\lambda}(w)$. 
This is entire. Therefore, $P_{\ell-1}(w) = 0$ on the zeros of $H^{(\ell)}_{\lambda}(w)$.
By Lemma \ref{hilbertnullstellensatz},
 there exists a polynomial $Q_{\ell}(w) $ such that 
\begin{equation*}
 \frac{P_{\ell-1}(w)}{H_{\lambda}^{(\ell)} (w)}=Q_{\ell}(w). 
\end{equation*}
Therefore, $F(w)/H_{\lambda}^{(1)} (w)\cdots H_{\lambda}^{(k)}(w)$ is a polynomial for $1 \leq k \leq N$. 
Taking $k = N$, we have that $F(w)/H_{\lambda}(w)$ is a polynomial of $w$, hence $f(z)/(h(z)-\lambda )$ is a polynomial of $e^{iz_j}$ by (\ref{variableschange}).
This implies that $\widehat{u}(n) $ is compactly supported.
We have thus completed the proof of Theorem \ref{rellich}.


\end{document}